\documentclass[10pt, a4paper]{article}
\usepackage[utf8]{inputenc}
\usepackage[russian]{babel}
\usepackage[T2A]{fontenc}
\usepackage{mathtext}
\usepackage{mathrsfs}
\usepackage{amssymb}
\usepackage{amsfonts}
\usepackage{tensor}
\usepackage{amscd}
\usepackage{amsmath, amsthm}
\usepackage{mathtools}
\usepackage{hyperref}
\usepackage[nice]{nicefrac}
\usepackage[symbol]{footmisc}
\usepackage{empheq}

\addto{\captionsrussian}{}

\newcommand*{\Af}{A_f^\text{\scriptsize{H}}}
\newcommand*\widebox[1]{\fbox{\hspace{1em}#1\hspace{1em}}}

\begin{document}

\begin{center}\textbf{TAKING THE LOGARITHM OF BINOMIAL TYPE SEQUENCES:\\LINEAR APPROACH}
\end{center}
\begin{center}
           \begin{footnotesize}
	DANIL KROTKOV\footnote[2]{Higher School of Economics, Faculty of Mathematics, e-mail: dikrotkov@hse.ru}
	\end{footnotesize}
\end{center}

\begin{footnotesize}
\noindent\textsc{Abstract}. In this paper we obtain the formal asymptotic expansion of the logarithms $\ln p_s(\alpha)$ of $p_s(\alpha)$, which are canonical continuations of polynomials of binomial type $p_n(\alpha)$. Our approach is based on linear methods which do not require the calculation of expansions $(p_s(\alpha)\alpha^{-s}-1)^k$, as opposed to the direct logarithmization.\\
\end{footnotesize}
~\\
In the previous paper (see $\cite{AfDf}$) for given series $f(x) \in x+x^2\mathbb{C}[[x]]$ and $D=\partial/\partial\alpha$ we suggest to define the complex powers of the operator $\alpha f'(D)^{-1}=A_f$ by the use of expression
$$
\Af \coloneqq\sum_{n=0}^\infty \binom{\text{\small{H}}-1}{n}\alpha^{\text{\scriptsize{H}}-n}q_n^D(\text{\small{H}})f'(D)^{-\text{\scriptsize{H}}} \eqno(1.1)
$$
where $q_n^t(s)$ are defined as coefficients of the expansion:
$$
\left(\frac{xf'(t)}{f(x+t)-f(t)}\right)^s=\sum_{n=0}^\infty \frac{q_n^{t}(s)x^n}{n!} \eqno(1.2)
$$
The following identities for powers may then be easily verified (here $D_f$ stands for $f(D)$)
$$
D_f^{\vphantom{s-1}} A_f^{s\vphantom{-1}}=sA_f^{s-1}+A_f^{s\vphantom{-1}}D_f^{\vphantom{s-1}}; ~~~~~~~~~~A_f^{s\vphantom{+1}} A_f^{h\vphantom{s+1}} =  A_f^{h\vphantom{s+1}}A_f^{s\vphantom{+1}} = A_f^{s+h\vphantom{+1}}
$$
Such an expression for $\Af$ may also be obtained by taking into account the following basic principle: Suppose there is an operator $T(\alpha, D)$, a family of series $g_t(\alpha)$, and $T(\alpha, D)\cdot g_t(\alpha)=\widetilde{g_t}(\alpha)$. Suppose there also exists an operator $\theta_g(\alpha, D)$, which does not depend on $t$, such that $\theta_g \cdot g_t(\alpha)=tg_t(\alpha)$. Consider then the expansion
$$
\frac{\widetilde{g_t}(\alpha)}{g_t(\alpha)}=\sum_{i\in I}\alpha^i \ell_i(t)
$$
(shifting, if necessary every $t$ to the right). Then if $\ell_i(t)$ are nice enough (for example, $\ell_i(t)\in\mathbb{C}[[t]]$, as in all cases considered below), one should speculate that the following equality of operators holds
$$
T(\alpha,D)\equiv\sum_{i\in I}\alpha^i \ell_i(\theta_g)
$$
It then follows that to obtain the expression for $\Af$, it is sufficient to know the image
$$
\Af \cdot e^{\alpha t}=p_\text{\scriptsize{H}}^{t\vphantom{()}}(\alpha)f'(t)^{-\text{\scriptsize{H}}}e^{\alpha t} \eqno(1.3)
$$
where
$$
p_\text{\scriptsize{H}}^{t\vphantom{()}}(\alpha)=\alpha \left(\frac{Df'(t)}{f(D+t)-f(t)}\right)^{\text{\scriptsize{H}}}\cdot\alpha^{{\text{\scriptsize{H}}-1}} \eqno(1.4)
$$
In this case we have $g_t(\alpha)=e^{\alpha t}$, $\theta_g(\alpha,D)=D$, and so the expansion (1.1) follows. But it should also be noticed that $\Af$ has another natural image (here $p_s(\alpha)=p_s^0(\alpha)$):
$$
A_f^{s}\cdot 1=p_s(\alpha) ~\Rightarrow~ \Af \cdot p_s(\alpha)=p_{s+\text{\scriptsize{H}}}(\alpha)
$$
(As a reminder, in case $s \in \mathbb{N}$ the sequence of polynomials $A_f^{s}\cdot 1=p_s(\alpha)$ is of binomial type $\cite{Dlt}$, and $\sum_{n=0}^\infty p_n(\alpha)f(x)^n/n!=\exp(\alpha x)$). Now since $A_fD_f\cdot p_s(\alpha)=sp_s(\alpha)$ (where $D_f$ stands for $f(D)$), one can use the basic principle again and consider the expansion:
$$
\frac{p_{s+\text{\scriptsize{H}}}(\alpha)}{p_{s}(\alpha)}=\sum_{n=0}^\infty \alpha^{\text{\scriptsize{H}}-n}P_n^{\text{\scriptsize{H}}}(s) \eqno(1.5)
$$
where $P_n^{\text{\scriptsize{H}}}(s)$ are polynomials of $s$ and $\text{\small{H}}$, what follows from the properties of expansions of $p_s(\alpha)$ and $p_{s+\text{\scriptsize{H}}}(\alpha)$. It then follows that:
$$
\Af=\sum_{n=0}^\infty \alpha^{\text{\scriptsize{H}}-n}P_n^{\text{\scriptsize{H}}}(A_fD_f) \eqno(1.6)
$$
Or, in other words, the following equality of operators holds
$$
\sum_{n=0}^\infty \binom{\text{\small{H}}-1}{n}\alpha^{\text{\scriptsize{H}}-n}q_n^D(\text{\small{H}})f'(D)^{-\text{\scriptsize{H}}}=\sum_{n=0}^\infty \alpha^{\text{\scriptsize{H}}-n}P_n^{\text{\scriptsize{H}}}(A_fD_f) \eqno(1.7)
$$
In the study of polynomials of binomial type and their continuations it is particularly interesting which properties has the fraction $p_{s+\text{\scriptsize{H}}}(\alpha)/p_{s}(\alpha)$ and its expansion.
We are interested in transforming the expansion of operator $T(\alpha,D)$ in powers of $\alpha$ and $D$ to its expansion in powers of $\alpha$ and $A_fD_f=\alpha f(D)f'(D)^{-1}$. Suppose this operator has the expansion of the form $T(\alpha, D)=\sum_{n\geqslant 0} g_n(\alpha)D^n$. Then we may perform the following symbolic calculation
\begin{align*}
\sum_{n=0}^\infty g_n(\alpha)D^n&=g_0(\alpha)+\sum_{n=0}^\infty g_{n+1}(\alpha)D^{n+1}\frac{f'(D)}{f(D)}\alpha^{-1}A_f D_f\\
&=g_0(\alpha)+\sum_{n=0}^\infty g_{n+1}(\alpha)\sum_{k=0}^\infty (-1)^{k}\alpha^{-1-k}\left(\frac{x^{n+1}f'(x)}{f(x)}\right)^{(k)}(D)A_f D_f
\end{align*}
We introduce the following operators, acting on $\mathbb{C}[[x]]$ (here $\mathrm{L}$ stands for the 0-derivative)
\begin{align*}
\mathrm{B}F(x)&\coloneqq \left(\alpha+\frac{d}{dx}\right)^{-1} f'(x)f(x)^{-1}(F(x)-F(0))=\left(\alpha+\frac{d}{dx}\right)^{-1}xf'(x)f(x)^{-1}\mathrm{L}\cdot F(x)\\
\mathrm{A}F(x)&\coloneqq F(0)
\end{align*}
It then follows that
$$
F(D)=(\mathrm{A}F)(D)+(\mathrm{B}F)(D)A_f D_f
$$
Continuing this procedure, we then obtain the following expression
$$
T\left(\alpha,\frac{\partial}{\partial\alpha}\right)=\left.\sum_{k=0}^\infty \alpha^{-k}\left[\left(1+\alpha^{-1}\frac{d}{dx}\right)^{-1}\frac{xf'(x)}{f(x)}\mathrm{L}~\right]^k \cdot ~T(\alpha,x)~\right|_{x=0} (A_fD_f)^k \eqno(1.8)
$$
It then follows that for the operator $\Af$ holds (here $\omega(x)$ is defined as an inverse series to $f(x)/f'(x)$, i.e. $f(\omega(x))/f'(\omega(x))=x$)
$$
\Af=\sum_{k=0}^\infty \alpha^{-k}\left.\left[\left(1+\alpha^{-1}\frac{d}{d\omega}\right)^{-1}\mathrm{L}~\right]^k \cdot~ p_{\text{\scriptsize{H}}}^{\omega(x)}(\alpha)f'(\omega(x))^{-\text{\scriptsize{H}}}~\right|_{x=0}(A_fD_f)^k
$$
Now act with both sides of this equality on $p_s(\alpha)$ and notice the geometric progression to obtain the identity
$$
\frac{p_{s+\text{\scriptsize{H}}}(\alpha)}{p_s(\alpha)}=\left.\left(1-\frac{s}{\alpha}\left(1+\alpha^{-1}\frac{d}{d\omega}\right)^{-1}\mathrm{L}\right)^{-1}\cdot ~p_{\text{\scriptsize{H}}}^{\omega(x)}(\alpha)f'(\omega(x))^{-\text{\scriptsize{H}}}~\right|_{x=0} \eqno(1.9)
$$
or equivalently
$$
\frac{p_{s+\text{\scriptsize{H}}}(\alpha)}{p_s(\alpha)}=\left.\left(1-s\alpha^{-1}\mathrm{L}+\alpha^{-1}\frac{d}{d\omega}\right)^{-1}\cdot~ \frac{p_{\text{\scriptsize{H}}+1}^{\omega(x)}(\alpha)}{\alpha}f'(\omega(x))^{-\text{\scriptsize{H}}}~\right|_{x=0} \eqno(1.10)
$$
Notice that in general, for $T(\alpha, D)$ the following holds true:
$$
\boxed{~~\frac{\alpha}{p_s(\alpha)}~T\cdot\frac{p_s(\alpha)}{\alpha}=\left.\left(1-s\alpha^{-1}\mathrm{L}+\alpha^{-1}\frac{d}{d\omega}\right)^{-1}e^{-\alpha\omega(x)}~T\cdot e^{\alpha\omega(x)} ~\right|_{x=0}~}
$$
Let us analyse this identity. First of all, since $\mathrm{L}=\left(1+x\frac{d}{dx}\right)^{-1}\frac{d}{dx}$, it may be rewritten as
$$
\frac{p_{s+\text{\scriptsize{H}}}(\alpha)}{p_s(\alpha)}=\left.\left(1+(x-s\alpha^{-1})\frac{d}{dx}+\alpha^{-1}\frac{d}{dx}x\frac{d}{d\omega}\right)^{-1}\frac{d}{dx}x\cdot~ \frac{p_{\text{\scriptsize{H}}+1}^{\omega(x)}(\alpha)}{\alpha}f'(\omega(x))^{-\text{\scriptsize{H}}}~\right|_{x=0} \eqno(1.11)
$$
where in brackets stands the hypergeometric-like operator (it reduces to hypergeometric one in case $f(x)=e^x-1 ~\Rightarrow~ \omega(x)=(f(x)/f'(x))^{inv}=-\ln(1-x)$). 
Consider now the formula (1.10). One may ask for whole representation of this expression, not only of its value at $x=0$. In fact,
\begin{align*}
&\left(1-s\alpha^{-1}\mathrm{L}+\alpha^{-1}\frac{d}{d\omega}\right)^{-1} \cdot~ g(\omega(x))=\\
&~~~~~=\frac{\alpha}{p_s(\alpha)}g(D)\cdot\frac{p_s(\alpha)}{\alpha}+\alpha e^{-\alpha\omega(x)}f(\omega(x))^s\int_{0}^{\omega(x)}e^{\alpha t}f(t)^{-s}\left(g(t)-\frac{\alpha}{p_s(\alpha)}g(D)\cdot\frac{p_s(\alpha)}{\alpha}\right)dt
\end{align*}
An interested reader may find the proof of this formula in Appendix A. Finally, one can expand (1.10) to obtain the following expression
\begin{align*}
&\frac{p_{s+\text{\scriptsize{H}}}(\alpha)}{p_s(\alpha)}=\sum_{n=0}^\infty \alpha^{-n}\left(s\mathrm{L}-\frac{d}{d\omega}\right)^n\left.\cdot~ \frac{p_{\text{\scriptsize{H}}+1}^{\omega(x)}(\alpha)}{\alpha}f'(\omega(x))^{-\text{\scriptsize{H}}}~\right|_{x=0}=\\
&~~~~=\sum_{n=0}^\infty\alpha^{\text{\scriptsize{H}}-n}\sum_{k=0}^n \left.\binom{\text{\small{H}}}{n-k}\left(s\mathrm{L}-\frac{d}{d\omega}\right)^k\cdot q_{n-k}^{\omega(x)}(1+\text{\small{H}})f'(\omega(x))^{-\text{\scriptsize{H}}}\right|_{x=0}
\end{align*}
where $q_k^t(s)$ are defined in (1.2). It then follows that polynomials $P_n^{\text{\scriptsize{H}}}(s)$ in the expression (1.5) when considered as polynomials in $s$ have the degree less or equal to $\leqslant n$, in spite of the fact that both the numerator and the denominator have the polynomials of degree $\leqslant 2n$ in $s$ as coefficients in their expansions. Now, one may take the operator $(1-s\alpha^{-1}\mathrm{L})^{-1}$ out of the brackets and use the following fact
$$
\tensor[_{x=0}]{\left|~(1-p\mathrm{L})^{-1\vphantom{\tfrac{d}{dx}}}~\cdot\right.}{} = \tensor[_{x=0}]{\left|~e^{p\tfrac{d}{dx}}~\cdot\right.}{}
$$
(since $(1-p\mathrm{L})^{-1}f(x)=(x-p)^{-1}(xf(x)-pf(p))$) to obtain the expression
$$
\frac{\alpha}{p_s(\alpha)}~T\cdot \frac{p_s(\alpha)}{\alpha}=\left.\left(1+\alpha^{-1}\frac{d}{d\omega}(1-s\alpha^{-1}\mathrm{L})^{-1}\right)^{-1}e^{-\alpha\omega(x)}~T\cdot e^{\alpha\omega(x)} ~\right|_{x=s\alpha^{-1}} \eqno(1.12)
$$
\noindent Multiply $s$ by $\alpha$:
$$
\left.\left[\frac{\alpha}{p_c(\alpha)}~T\cdot \frac{p_c(\alpha)}{\alpha}\right]\right|_{c=s\alpha}=\left.\left(1+\alpha^{-1}\frac{d}{d\omega}\left(1+(x-s)\frac{d}{dx}\right)^{-1}\frac{d}{dx}x\right)^{-1}e^{-\alpha\omega(x)}~T\cdot e^{\alpha\omega(x)} ~\right|_{x=s}
$$
Now make the change of variables $x \to x+s$ and rewrite $f(x+s)$ as $e^{\strut x\tfrac{\partial}{\partial s}}f(s)e^{-\strut x\tfrac{\partial}{\partial s}}$, to obtain
\begin{align*}
\left.\left[\frac{\alpha}{p_c(\alpha)}~T\cdot \frac{p_c(\alpha)}{\alpha}\right]\right|_{c=s\alpha}&=\left.\left(1+\frac{\alpha^{-1}}{\omega'(x+s)}\frac{d}{dx}\mathrm{L}~(x+s)\right)^{-1}e^{-\alpha\omega(x+s)}~T\cdot e^{\alpha\omega(x+s)} ~\right|_{x=0}=\\
&=\left.\left(1+\frac{\alpha^{-1}}{\omega'(x+s)}\frac{d}{dx}\mathrm{L}~(x+s)\right)^{-1} e^{\strut x\tfrac{\partial}{\partial s}} ~\right|_{x=0} e^{-\alpha\omega(s)}~T\cdot e^{\alpha\omega(s)}=\\
&=\frac{1}{x+s}\left.\left(1+\frac{\alpha^{-1}(x+s)}{\omega'(x+s)}\frac{d}{dx}\mathrm{L}\right)^{-1} e^{\strut x\tfrac{\partial}{\partial s}} ~\right|_{x=0} s e^{-\alpha\omega(s)}~T\cdot e^{\alpha\omega(s)}=\\
&=s^{-1}\left.\left(1+\frac{\alpha^{-1}(x+s)}{\omega'(x+s)}\frac{d}{dx}\mathrm{L}\right)^{-1} e^{\strut x\tfrac{\partial}{\partial s}} ~\right|_{x=0} s e^{-\alpha\omega(s)}~T\cdot e^{\alpha\omega(s)}\\
&=s^{-1}\left.\left(1+\frac{\alpha^{-1}s}{\omega'(s)}e^{-\strut x\tfrac{\partial}{\partial s}}\frac{d}{dx}\mathrm{L}~e^{\strut x\tfrac{\partial}{\partial s}}\right)^{-1} ~\right|_{x=0} s e^{-\alpha\omega(s)}~T\cdot e^{\alpha\omega(s)}
\end{align*}
Now introduce the operator
$$
M_\alpha\left(s,\frac{\partial}{\partial s}\right)=\left.\left(1+\alpha^{-1}\frac{s}{\omega'(s)}e^{-\strut x\tfrac{\partial}{\partial s}}\frac{d}{dx}\mathrm{L}~e^{\strut x\tfrac{\partial}{\partial s}}\right)^{-1} ~\right|_{x=0}
$$
Hence, the following holds:
$$
\left.\left[\frac{\alpha}{p_c(\alpha)}~T\left(\alpha,\frac{\partial}{\partial\alpha}\right)\cdot \frac{p_c(\alpha)}{\alpha}\right]\right|_{c=s\alpha}=s^{-1}M_\alpha\left(s,\frac{\partial}{\partial s}\right)se^{-\alpha\omega(s)}~T\left(\alpha,\frac{\partial}{\partial\alpha}\right)\cdot e^{\alpha\omega(s)}\eqno(1.13)
$$
From now on we are interested in this operator. Expanding it in powers of $\alpha$, we have
$$
M_\alpha\left(s,\frac{\partial}{\partial s}\right)=\sum_{n=0}^\infty (-\alpha)^{-n}~T_n\left(s,\frac{\partial}{\partial s}\right)
$$
where
$$
T_n\left(s,\frac{\partial}{\partial s}\right)f(s)=\left.\left(\frac{s}{\omega'(s)}e^{-\strut x\tfrac{\partial}{\partial s}}\frac{d}{dx}\mathrm{L}~ e^{\strut x\tfrac{\partial}{\partial s}} \right)^n\right|_{x=0}f(s)=\tensor[_{x=0}]{\left|~e^{\strut s\tfrac{d}{dx}} \left(x\frac{d}{d\omega}(1-s\mathrm{L})^{-1}\mathrm{L}\right)^n f(x)\right.}{}
$$
\noindent\underline{\textbf{Theorem}}. Consider the structure of noncommutative polynomials in variables $\sigma, E, D$ of the following form:
$$
\alpha_0 E+\alpha_1 ED+...+\alpha_k ED^k
$$
where coefficients $\alpha_i$ depend on $\sigma, E, D$ somehow (note that we are interested only in polynomials of this form, i.e. without the \emph{true} free term). Consider the following transformation (which does not affect $\alpha_i$):
$$
ED^i \xrightarrow{~~\nu~~} \frac{1}{(i+1)(i+2)}\sigma D^{i+2}E-\frac{1}{i+1}\sigma DED^{i+1}+\frac{1}{i+2}\sigma ED^{i+2}
$$
It means that the image of the polynomial is of the form $\widetilde{\alpha}_0E+\widetilde{\alpha}_1ED+...+\widetilde{\alpha}_{k+2}ED^{k+2}$, where $\widetilde{\alpha}_i$ are determined by conditions:
\begin{align*}
&\widetilde{\alpha}_0=\frac{1}{1\cdot2}\alpha_0\sigma D^2+\frac{1}{2\cdot3}\alpha_1\sigma D^3+...+\frac{1}{(k+1)(k+2)}\alpha_{k}\sigma D^{k+2}\\
&\widetilde{\alpha}_1=-\alpha_0\sigma D\\
&\widetilde{\alpha}_{i}=\frac{1}{i}(\alpha_{i-2}\sigma-\alpha_{i-1}\sigma D) ~~~~(1<i\leqslant n+1)\\
&\widetilde{\alpha}_{k+2}=\frac{1}{k+2}\alpha_{k}\sigma
\end{align*}
Now set the initial polynomial $P_0=E$. Suppose $P_n=\nu^n\cdot P_0$, i.e. $P_1=\frac{1}{2}\sigma D^2 E - \sigma DED+\frac{1}{2}\sigma ED^2$ $~~~(=\frac{1}{2}\sigma [D,[D,E]])$ and so on. $\vphantom{\widetilde{T}_1=\tfrac{1}{2}\tfrac{s}{\omega'(s)}\tfrac{\partial^2}{\partial s^2}}$ Consider the term $\alpha_0$ of each of $P_n$ and set $D =\frac{\partial}{\partial s}$, $\sigma=\frac{s}{\omega'(s)}$. In other words $P_0 ~\to~ \widetilde{T}_0=1$, $P_1=\frac{1}{2}\sigma D^2 ~\to~ \widetilde{T}_1=\tfrac{1}{2}\tfrac{s}{\omega'(s)}\tfrac{\partial^2}{\partial s^2}$, ... . Then the following holds
$$
\widetilde{T}_n=T_n\left(s,\frac{\partial}{\partial s}\right)
$$
In order to prove this theorem consider the following proposition.\\
~\\
\underline{\textbf{Proposition}}. Suppose $\mathrm{L}f(x)=x^{-1}(f(x)-f(0))$. Then
$$
\tensor[_{x=0}]{\left|~e^{\strut p\tfrac{d}{dx}} \frac{d^n}{dp^n}(1-p\mathrm{L})^{-1}\mathrm{L}\right.}{}=\tensor[_{x=0}]{\left|\frac{1}{n+1}\left[\frac{d^{n+1}}{dp^{n+1}}e^{\strut p\tfrac{d}{dx}}-e^{\strut p\tfrac{d}{dx}}\frac{d^{n+1}}{dp^{n+1}}\right]\right.}{}
$$
\emph{Proof}: We have
$$
(1-p\mathrm{L})^{-1}f(x)=\frac{xf(x)-pf(p)}{x-p}
$$
Hence
$$
(1-p\mathrm{L})^{-1}=\int_{0}^{1}t^{\strut x\tfrac{d}{dx}}e^{p(1-t)\strut\tfrac{d}{dx}}\frac{d}{dx}xdt ~\Rightarrow~ (1-p\mathrm{L})^{-1}\mathrm{L}=\int_{0}^{1}t^{\strut x\tfrac{d}{dx}}e^{p(1-t)\strut\tfrac{d}{dx}}\frac{d}{dx}x\mathrm{L}~dt \Rightarrow
$$
$$
\Rightarrow (1-p\mathrm{L})^{-1}\mathrm{L}=\int_{0}^{1}t^{\strut x\tfrac{d}{dx}}e^{p(1-t)\strut\tfrac{d}{dx}}\frac{d}{dx}dt
$$
(here annihilation of $x\mathrm{L}\neq 1$ happens because $\frac{d}{dx}x\frac{f(x)-f(0)}{x}=\frac{d}{dx}f(x)$). The latter means that
\begin{align*}
&\tensor[_{x=0}]{\left|~e^{\strut p\tfrac{d}{dx}} \frac{d^n}{dp^n}(1-p\mathrm{L})^{-1}\mathrm{L}\right.}{}=\tensor[_{x=0}]{\left|~e^{\strut p\tfrac{d}{dx}} \frac{d^n}{dp^n}\int_{0}^{1}t^{\strut x\tfrac{d}{dx}}e^{p(1-t)\strut\tfrac{d}{dx}}\frac{d}{dx}dt\right.}{}=\\
&~~~=\tensor[_{x=0}]{\left|~e^{\strut p\tfrac{d}{dx}} \int_{0}^{1}t^{\strut x\tfrac{d}{dx}}e^{p(1-t)\strut\tfrac{d}{dx}}\left(\frac{d}{dp}+(1-t)\frac{d}{dx}\right)^n\frac{d}{dx}dt\right.}{}\\
&~~~=\tensor[_{x=0}]{\left|~ \int_{0}^{1}e^{\strut pt\tfrac{d}{dx}}e^{p(1-t)\strut\tfrac{d}{dx}}\left(\frac{d}{dp}+(1-t)\frac{d}{dx}\right)^n\frac{d}{dx}dt\right.}{}=\\
&~~~=\tensor[_{x=0}]{\left|~ \int_{0}^{1}e^{p\strut\tfrac{d}{dx}}\left(\frac{d}{dp}+(1-t)\frac{d}{dx}\right)^n\frac{d}{dx}dt\right.}{}=\tensor[_{x=0}]{\left|~\sum_{k=0}^{n}\binom{n}{k}\frac{1}{k+1}e^{p\strut\tfrac{d}{dx}}\frac{d^{n-k}}{dp^{n-k}}\frac{d^{k+1}}{dx^{k+1}}\right.}{}=\\
\tag*{\qed}&~~~=\frac{1}{n+1}\tensor[_{x=0}]{\left|~e^{p\strut\tfrac{d}{dx}}\left(\left(\frac{d}{dx}+\frac{d}{dp}\right)^{n+1}-\frac{d^{n+1}}{dp^{n+1}}\right)\right.}{}=\tensor[_{x=0}]{\left|\frac{1}{n+1}\left[\frac{d^{n+1}}{dp^{n+1}}e^{\strut p\tfrac{d}{dx}}-e^{\strut p\tfrac{d}{dx}}\frac{d^{n+1}}{dp^{n+1}}\right]\right.}{}
\end{align*}
Then suppose that there exist some $\alpha_i^n(x,y)$, such that
$$
\tensor[_{x=0}]{\left|~e^{s\strut\tfrac{d}{dx}}\left(\frac{x}{\omega'(x)}\frac{d}{dx}(1-s\mathrm{L})^{-1}\mathrm{L}\right)^n\right.}{}=\tensor[_{x=0}]{\left|~\sum_i \alpha_i^n\left(s,\frac{\partial}{\partial s}\right) e^{s\strut\tfrac{d}{dx}}\frac{\partial^i}{\partial s^i}\right.}{}
$$
As in case $n=1$:
\begin{align*}
&\tensor[_{x=0}]{\left|~e^{s\strut\tfrac{d}{dx}}\frac{x}{\omega'(x)}\frac{d}{dx}(1-s\mathrm{L})^{-1}\mathrm{L}\right.}{}=\tensor[_{x=0}]{\left|~\frac{s}{\omega'(s)}e^{s\strut\tfrac{d}{dx}}\frac{d}{dx}(1-s\mathrm{L})^{-1}\mathrm{L}\right.}{}=\\
&~~~=\tensor[_{x=0}]{\left|~\frac{s}{\omega'(s)}\left(\frac{\partial}{\partial s}e^{s\strut\tfrac{d}{dx}}-e^{s\strut\tfrac{d}{dx}}\frac{\partial}{\partial s}\right)(1-s\mathrm{L})^{-1}\mathrm{L}\right.}{}=\\
&~~~=\tensor[_{x=0}]{\left|~\frac{s}{\omega'(s)}\left(\frac{1}{2}\frac{\partial^2}{\partial s^2}e^{s\strut\tfrac{d}{dx}}-\frac{\partial}{\partial s}e^{s\strut\tfrac{d}{dx}}\frac{\partial}{\partial s}+\frac{1}{2}e^{s\strut\tfrac{d}{dx}}\frac{\partial^2}{\partial s^2}\right)\right.}{}
\end{align*}
Then, in general, we have
\begin{align*}
&\tensor[_{x=0}]{\left|~e^{s\strut\tfrac{d}{dx}}\left(\frac{x}{\omega'(x)}\frac{d}{dx}(1-s\mathrm{L})^{-1}\mathrm{L}\right)^{n+1}\right.}{}=\tensor[_{x=0}]{\left|~\sum_i \alpha_i^n\left(s,\frac{\partial}{\partial s}\right) e^{s\strut\tfrac{d}{dx}}\frac{\partial^i}{\partial s^i}\frac{x}{\omega'(x)}\frac{d}{dx}(1-s\mathrm{L})^{-1}\mathrm{L}\right.}{}=\\
&=\tensor[_{x=0}]{\left|~\sum_i \alpha_i^n\left(s,\frac{\partial}{\partial s}\right)\frac{s}{\omega'(s)}e^{s\strut\tfrac{d}{dx}}\frac{\partial^i}{\partial s^i}\frac{d}{dx}(1-s\mathrm{L})^{-1}\mathrm{L}\right.}{}=\\
&=\tensor[_{x=0}]{\left|~\sum_i \alpha_i^n\left(s,\frac{\partial}{\partial s}\right)\frac{s}{\omega'(s)}\left(\frac{\partial}{\partial s}e^{s\strut\tfrac{d}{dx}}-e^{s\strut\tfrac{d}{dx}}\frac{\partial}{\partial s}\right)\frac{\partial^i}{\partial s^i}(1-s\mathrm{L})^{-1}\mathrm{L}\right.}{}=\tensor[_{x=0}]{\left|~\sum_i \alpha_i^n\left(s,\frac{\partial}{\partial s}\right)\right.}{}\times\\
&~~~~~~~~~~~~~~\times\frac{s}{\omega'(s)}\left(\frac{1}{(i+1)(i+2)}\frac{\partial^{i+2}}{\partial s^{i+2}}e^{s\strut\tfrac{d}{dx}}-\frac{1}{i+1}\frac{\partial}{\partial s}e^{s\strut\tfrac{d}{dx}}\frac{\partial^{i+1}}{\partial s^{i+1}}+\frac{1}{i+2}e^{s\strut\tfrac{d}{dx}}\frac{\partial^{i+2}}{\partial s^{i+2}}\right)=\\
&=\tensor[_{x=0}]{\left|~\sum_i \alpha_i^{n+1}\left(s,\frac{\partial}{\partial s}\right) e^{s\strut\tfrac{d}{dx}}\frac{\partial^i}{\partial s^i}\right.}{}
\end{align*}
It is left to say that with action on $f(x)$ the only term that does not vanish is $\alpha_0^{n}\left(s,\frac{\partial}{\partial s}\right)\cdot f(s)$. Now, since the coefficients do indeed transform like:
\begin{align*}
&\alpha_0^{n+1}=\sum_i\frac{1}{(i+1)(i+2)}\alpha_i^n\frac{s}{\omega'(s)}\frac{\partial^{i+2}}{\partial s^{i+2}}\\
&\alpha_1^{n+1}=-\alpha_0^n\frac{s}{\omega'(s)}\frac{\partial}{\partial s}\\
&\alpha_{i}^{n+1}=\frac{1}{i}\alpha_{i-2}^n\frac{s}{\omega'(s)}-\frac{1}{i}\alpha_{i-1}^n\frac{s}{\omega'(s)} \frac{\partial}{\partial s} ~~~~(1<i)
\end{align*}
we obtain the desired result. \hfill\ensuremath{\blacksquare}\\
~\\
The latter approach gives us a simple method (at least theoretically) to evaluate $T_i$. So $T_0=1$, $T_1(s, \frac{\partial}{\partial s})=\frac{1}{2}\frac{s}{\omega'(s)}\frac{\partial^2}{\partial s^2}$, $T_2(s, \frac{\partial}{\partial s})=\frac{1}{4}\frac{s}{\omega'(s)}\frac{\partial^2}{\partial s^2}\frac{s}{\omega'(s)}\frac{\partial^2}{\partial s^2}-\frac{1}{6}\frac{s}{\omega'(s)}\frac{\partial}{\partial s}\frac{s}{\omega'(s)}\frac{\partial^3}{\partial s^3}+\frac{1}{24}\frac{s^2}{\omega'(s)^2}\frac{\partial^4}{\partial s^4}$. In other\\ words, $T_n(s, \frac{\partial}{\partial s})$ may be obtained by successive multiplication of the following $(2k-1)\times(2k+1)$-matrices:
$$
\begin{pmatrix}
\frac{1}{2}\sigma D^2 &-\sigma D & \frac{1}{2}\sigma &0&0&0&\hdots&0\\
~\\
\frac{1}{6}\sigma D^3 &0&-\frac{1}{2}\sigma D&\frac{1}{3}\sigma&0&0&\hdots&0\\
~\\
\frac{1}{12}\sigma D^4 &0&0&-\frac{1}{3}\sigma D&\frac{1}{4}\sigma&0&\hdots&0\\
~\\
\frac{1}{20}\sigma D^5 &0&0&0&-\frac{1}{4}\sigma D&\frac{1}{5}\sigma&\hdots&0\\
\vdots&\vdots&\vdots&\vdots&\vdots&\vdots&\ddots\\
\frac{1}{(2k-1)2k}\sigma D^{2k}&0&0&0&0&0&{}&\frac{1}{2k}\sigma
\end{pmatrix}
$$
and by taking the first element of the resulting row. 
It is left to apply the obtained methods to the operator $\alpha^{-1}\Af \alpha$.
\begin{align*}
&\frac{p_{\alpha s+\text{\scriptsize{H}}}^{\vphantom{\omega}}(\alpha)}{p_{\alpha s}(\alpha)}=s^{-1}M_{\alpha}\left(s,\frac{\partial}{\partial s}\right)s~\frac{p_{\text{\scriptsize{H}}+1}^{\omega(s)}(\alpha)}{\alpha}f'(\omega(s))^{-\text{\scriptsize{H}}} ~\Rightarrow\\
&~~~\Rightarrow~\frac{p_{\alpha s+\text{\scriptsize{H}}}^{\vphantom{\omega}}(\alpha)}{p_{\alpha s}(\alpha)}=s^{-1}\sum_{n=0}^\infty (-\alpha)^{-n}~T_n\left(s,\frac{\partial}{\partial s}\right)s\sum_{n=0}^\infty \alpha^{\text{\scriptsize{H}}-n}\binom{\text{\small{H}}}{n}q_n^{\omega(s)}(1+\text{\small{H}})f'(\omega(s))^{-\text{\scriptsize{H}}}=\\
&~~~~~~=\sum_{n=0}^\infty \alpha^{\text{\scriptsize{H}}-n}\sum_{k=0}^n \binom{\text{\small{H}}}{n-k}(-1)^ks^{-1}T_k\left(s,\frac{\partial}{\partial s}\right)sq_{n-k}^{\omega(s)}(1+\text{\small{H}})f'(\omega(s))^{-\text{\scriptsize{H}}}\approx\\
&~~~~~~\approx\alpha^{\text{\scriptsize{H}}}f'(\omega(s))^{-\text{\scriptsize{H}}}+\alpha^{\text{\scriptsize{H}}-1}\left(\text{\small{H}}q_1^{\omega(s)}(1+\text{\small{H}})f'(\omega(s))^{-\text{\scriptsize{H}}}-s^{-1}\frac{1}{2}\frac{s}{\omega'(s)}\frac{\partial^2}{\partial s^2}sf'(\omega(s))^{-\text{\scriptsize{H}}}\right)
\end{align*}
By definition (1.2) of $q_k^t(s)$, $q_0^t(s)=1$ and
$$
q_1^t(s)=-\frac{s}{2}\frac{f''(t)}{f'(t)} ~\Rightarrow~ q_1^{\omega(s)}(1+\text{\small{H}})=\frac{1}{2}(1+\text{\small{H}})\frac{1-\omega'(s)}{s\omega'(s)}
$$
Hence
$$
\frac{p_{\alpha s+\text{\scriptsize{H}}}^{\vphantom{\omega}}(\alpha)}{p_{\alpha s}(\alpha)}\approx \alpha^{\text{\scriptsize{H}}}f'(\omega(s))^{-\text{\scriptsize{H}}}+\alpha^{\text{\scriptsize{H}}-1}\left(\frac{\text{\small{H}}^2}{2}\frac{1-\omega'(s)}{s}f'(\omega(s))^{-\text{\scriptsize{H}}}+\frac{\text{\small{H}}}{2}\frac{\omega''(s)}{\omega'(s)}f'(\omega(s))^{-\text{\scriptsize{H}}}\right)
$$
So, the latter implies the following limit formulae (here $\dot p_s(\alpha)=\frac{\partial}{\partial s} p_s(\alpha)$).
\begin{center}
\textbf{\underline{First Limit Formula}}:
\begin{align*}
\lim_{s\to\infty} \frac{\dot p_{s}^{\vphantom{\omega}}(s\alpha)}{p_{s}(s\alpha)}-\ln s=\ln\alpha-\ln f'(\omega(\alpha^{-1}))~~~~~~~~\lim_{s\to\infty} \frac{p_{s+\text{\scriptsize{H}}}^{\vphantom{\omega}}(s\alpha)}{p_{s}(s\alpha)}s^{-\text{\scriptsize{H}}}=\alpha^{\text{\scriptsize{H}}}f'(\omega(\alpha^{-1}))^{-\text{\scriptsize{H}}}
\end{align*}
\textbf{\underline{Second Limit Formula}}:
\begin{align*}
\lim_{s\to\infty} \ln p_s(\alpha s)-s\ln (\alpha s)&+s\alpha\int_{0}^{\alpha^{-1}}\ln f'(\omega(t))dt=\frac{1}{2}\ln\omega'(\alpha^{-1})\\
\lim_{s\to\infty} \frac{p_{s+\text{\scriptsize{H}}}^{\vphantom{\omega}}(s\alpha)}{p_{s}(s\alpha)}s^{1-\text{\scriptsize{H}}}-s\alpha^{\text{\scriptsize{H}}}f'(\omega(\alpha^{-1}))^{-\text{\scriptsize{H}}}&=\frac{\text{\small{H}}}{2}\alpha^{\text{\scriptsize{H}}}f'(\omega(\alpha^{-1}))^{-\text{\scriptsize{H}}}\left(\text{\small{H}}(1-\omega'(\alpha^{-1}))+\alpha^{-1}\frac{\omega''(\alpha^{-1})}{\omega'(\alpha^{-1})}\right)
\end{align*}
\end{center}
In particular case of the logarithm it is more useful to consider an operator $Df'(D)/f(D)$ instead of $\alpha^{-1}\ln A_f ~\alpha$, since there is an equality
$$
\frac{\alpha}{p_s(\alpha)}\frac{Df'(D)}{f(D)}\cdot \frac{p_s(\alpha)}{\alpha}=\frac{\alpha}{s}\frac{p_s'(\alpha)}{p_s(\alpha)}
$$
Hence
$$
\frac{p_{s\alpha}'(s)}{p_{s\alpha}(s)}=M_s\left(\alpha, \frac{\partial}{\partial \alpha}\right)\cdot \omega(\alpha)
$$
And so
$$
-\frac{s}{\alpha^2}\frac{p_s'}{p_s}\left(\frac{s}{\alpha}\right)=\sum_{n=0}^\infty (-s)^{1-n} \alpha^{n-2} T_n\left(\alpha, \frac{\partial}{\partial \alpha}\right)\cdot \omega(\alpha)
$$
Thus, integrating by $\alpha$, one can obtain the following expansion (generalized Stirling's formula):
\begin{subequations}
\begin{empheq}[box=\widebox]{align*}
&~\\
&\ln p_s(s\alpha^{-1})\sim s\ln (s\alpha^{-1})-s\alpha^{-1}\int_{0}^\alpha\ln f'(\omega(t))dt+\frac{1}{2}\ln\omega'(\alpha)+\\
&~~~~~+\frac{1}{24s}\left(2\frac{\omega'(\alpha)-1}{\omega'(\alpha)}+4\frac{\alpha^2\omega''(\alpha)^2}{\omega'(\alpha)^3}-2\frac{\alpha\omega''(\alpha)}{\omega'(\alpha)^2}-3\frac{\alpha^2 \omega^{(3)}(\alpha)}{\omega'(\alpha)^2}\right)-\frac{1}{48s^2}\frac{\alpha^{3}}{\omega'(\alpha)}\left(\frac{\alpha}{\omega'(\alpha)}\right)^{(4)}+...\\
\end{empheq}
\end{subequations}
which is equivalent to the identity
$$
\frac{1}{s}\ln (\alpha^{-s}p_s(\alpha)) = \ln\left(1-s\alpha^{-1}\mathrm{L}+\alpha^{-1}\frac{d}{d\omega}\right)\cdot\left.\frac{\omega(x)}{x}~\right|_{x=0}\eqno(1.14)
$$
or
$$
\frac{1}{s}\ln p_s(\alpha)=e^{\tfrac{\partial}{\partial\alpha}\left(\tfrac{d}{dx}-s\tfrac{xf'(x)}{f(x)}\mathrm{L}\right)}\cdot\left.\frac{xf'(x)}{f(x)}~\right|_{x=0} \cdot ~\ln \alpha \eqno(1.15)
$$
It is interesting that each term $g_n(\alpha)$ of this expansion except the first two ones is invariant under the transformation $f(x) \to \widetilde{f}(x)=f(x)e^{-Ax};~~\alpha \to \alpha(1+A\alpha)^{-1}$, since for $f \rightsquigarrow p_s$, $\widetilde{f} \rightsquigarrow \widetilde{p}_s$ holds:
$$
\alpha \widetilde{p}_s(s\alpha^{-1})=\frac{\alpha}{1+A\alpha} p_s(s\alpha^{-1}+As)
$$
(Here we use the expression $f(x) \rightsquigarrow p_s(\alpha)$ for $f(x) \in x+x^2\mathbb{C}[[x]]$ to denote the formal series $p_s(\alpha) \in \alpha^s+\alpha^{s-1}\mathbb{C}[[\alpha^{-1}]]$, which is determined by condition $p_s(\alpha)=\alpha \left(D/f(D)\right)^s \alpha^{s-1}$)
and at the same time
$$
\widetilde{\omega}(x)=(\mathfrak{T}(f(x)e^{-Ax}))^{inv}=\omega\left(\frac{x}{1+Ax}\right)
$$
It should also be mentioned that in case of polynomials $xe^{\mathrm{Ei}(x)-\ln|x|-\gamma}=\mathfrak{T}^{-1}xe^{-x}\rightsquigarrow \nu_k(\alpha)$ and their continuations $\nu_s(\alpha)$ the first terms of this formal asymptotic expansion may be written as
$$
\ln ((\alpha s^{-1})^s\nu_s(s \alpha^{-1})) ~\sim~ -s\sum_{n=1}^\infty \frac{\alpha^n}{n}\frac{(n+1)^{n-1}}{n!}+\frac{1}{2}\sum_{n=1}^\infty \frac{\alpha^n}{n}\sum_{k=0}^{n} \frac{n^k}{k!}+...
$$
\begin{center}
\textbf{General comment}
\end{center}
The next natural step is to ask how to take the logarithms of arbitrary Sheffer sequences $\cite{Sh}$. It is easy to verify that in this case the general arguments, we use in this paper, remain the same except minor changes. To be more specific, consider the series $\ell(x)\in1+x\mathbb{C}[[x]]$ and the corresponding series of the form $\alpha^s+\alpha^{s-1}\mathbb{C}[[\alpha^{-1}]]$:
$$
\tau_s^{f,\ell}(\alpha)\coloneqq\ell(D)A_f^{s}\cdot 1
$$
Which are polynomials in case $s \in \mathbb{N}$, since
$$
\sum_{n=0}^\infty \frac{\tau_n^{f,\ell}(\alpha)}{n!}x^n=\ell(\varphi(x))e^{\alpha \varphi(x)}
$$
(here $\varphi(x)=f(x)^{inv}$). Then one may consider the operator $\theta_{f,\ell}$:
$$
\theta_{f,\ell}\coloneqq\ell(D)A_fD_f\ell(D)^{-1} ~~\Rightarrow~~ \theta_{f,\ell} \tau_s^{f,\ell}(\alpha)=s\tau_s^{f,\ell}(\alpha)
$$
Consider again $\omega(x)=(\mathfrak{T}f(x))^{inv}$ and introduce the following two auxiliary operators on $\mathbb{C}[[x]]$:
\begin{align*}
&\widetilde{\mathrm{A}}F(x)\coloneqq F(0)\\
&\widetilde{\mathrm{B}}F(x)\coloneqq \alpha^{-1}\ell(\omega(x))^{-1}\left(1+\alpha^{-1}\frac{d}{d\omega}\right)^{-1} \ell(\omega(x))\mathrm{L}~F(x)
\end{align*}
It then follows that
$$
F(\mathfrak{T}f(D))=(\widetilde{\mathrm{A}}F)(\mathfrak{T}f(D))+(\widetilde{\mathrm{B}}F)(\mathfrak{T}f(D))\theta_{f,\ell}
$$
And hence for the operator $T(\alpha, D)=\sum \alpha^{i} q_i(D)$ and the corresponding series $T(\alpha, \omega)=\sum \alpha^{i} q_i(\omega)$ the following identity holds true:
$$
\frac{\alpha}{\tau_s^{f,\ell}(\alpha)}T(\alpha,D)\tau_{s-1}^{f,\ell}(\alpha)=\left.\left(1+\alpha^{-1}\frac{d}{d\omega}-s\alpha^{-1}\ell(\omega)\mathrm{L}~\ell(\omega)^{-1}\right)^{-1}\cdot ~T(\alpha,\omega(x))\ell(\omega(x))f'(\omega(x))~\right|_{x=0}
$$
One may consider again the following expansion
$$
\left.\frac{\partial}{\partial x}\ln\tau_s^{f,\ell}(x)~\right|_{x=s\alpha^{-1}}=\left.\sum_{n=0}^\infty (-\alpha)^ns^{-n}\left(\ell(\omega)^{-1}x\frac{d}{d\omega}\ell(\omega)(1-\alpha\mathrm{L})^{-1}\mathrm{L}\right)^n \cdot ~\omega(x)~\right|_{x=\alpha}
$$
and notice that the operators
$$
T_n^\ell\left(\alpha,\frac{\partial}{\partial \alpha}\right)\cdot g(\alpha)\coloneqq\tensor[_{x=0}]{\left|~e^{\alpha\frac{d}{dx}}\left(\ell(\omega)^{-1}x\frac{d}{d\omega}\ell(\omega)(1-\alpha\mathrm{L})^{-1}\mathrm{L}\right)^n\cdot~ g(x)~\right.}{}
$$
may be obtained by the use of structure of noncommutative polynomials in variables $\sigma,E,D,\lambda$ (where $\sigma=\frac{\alpha}{\omega'(\alpha)}, E=e^{\alpha \frac{d}{dx}}, D=\frac{\partial}{\partial \alpha}, \lambda=\ell(\omega(\alpha))$) and the transformation $\bar{\nu}$ acting on it:
$$
ED^i \xrightarrow{~~\bar{\nu}~~} \left(\frac{1}{i+1}\sigma \lambda^{-1}D\lambda D^{i+1}-\frac{1}{i+2}\sigma D^{i+2} \right)E-\frac{1}{i+1}\sigma \lambda^{-1}D\lambda ED^{i+1}+\frac{1}{i+2}\sigma ED^{i+2}
$$
So, again $T_n^{\ell}(\alpha,\frac{\partial}{\partial\alpha})$ is the first element of the polynomial $\bar{\nu}^n\cdot E$. Using this principle one may derive, as an example, the following representation of the logarithms of Bernoulli polynomials and their continuations:
$$
\frac{1}{s}\ln(B_s(\alpha)\alpha^{-s})=\left.\ln\left(1+\alpha^{-1}\frac{d}{dx}-s\alpha^{-1}\frac{x}{e^x-1}\mathrm{L} ~\frac{e^x-1}{x}\right)\cdot \frac{x}{e^x-1}\right|_{x=0}
$$
If one wants to make even more general statements, then one may consider an operator $\Gamma=$ $=\sum_{n=0}^\infty \alpha^{q-n}\mathfrak{g}_n(D)$ (where $\mathfrak{g}_0(0)\neq 0$, $\mathfrak{g}_n(t)\in \mathbb{C}[[t]]$). Then there exist some $g_i(t) \in \mathbb{C}[[t]]$ such that
$$
\Gamma^{-1}=\left(\sum_{n=0}^\infty \alpha^{q-n}\mathfrak{g}_n(D)\right)^{-1}=\sum_{n=0}^\infty g_n(D)\alpha^{-q-n} ~~\xRightarrow{~def~}~~ \overline{\Gamma^{-1}}\coloneqq\sum_{n=0}^\infty \alpha^{-q-n}g_n(-D)
$$
Now if $p_s(\alpha)$ is such a sequence of series that
$$
\Gamma D\cdot p_s(\alpha)=sp_s(\alpha)
$$
then the following statement holds true:
$$
\frac{1}{p_s(\alpha)}T(\alpha, D) \cdot p_s(\alpha)=\left.\left(1-s~e^{\frac{\partial^2}{\partial \alpha\partial x}}e^{\alpha x}~\overline{\Gamma^{-1}}~e^{-\alpha x}e^{-\frac{\partial^2}{\partial \alpha\partial x}}~\mathrm{L}\right)^{-1} e^{-\alpha x}~T(\alpha,D)\cdot e^{\alpha x}~\right|_{x=0}
$$
\begin{center}
\textbf{Conclusion}
\end{center}

Such an approach, introduced in this paper allows us to formulate theorems, which may not be seen at first glance, as for example:
$$
\lim_{n \to \infty}\frac{p_n'(n\alpha)}{p_n(n\alpha)}=\omega(\alpha^{-1})
$$
Nevertheless, the introduced expansions are not as well understood as in classical case $f(x)=e^x-1 \rightsquigarrow p_s(\alpha)=(\alpha)_s$, and so there is still a huge amount of work to be done.

\newpage

\begin{center}
\textbf{Appendix A}
\end{center}

\hrule
~\\
~\\

\noindent Consider the expansion
$$
\left(\frac{x}{f(x)}\right)^s=\sum_{n=0}^\infty \frac{q_n(s)x^n}{n!}
$$
\textbf{Proposition A.1}.
$$
\left(\frac{\omega}{f(\omega)}\right)^s\left(s\mathrm{L}-\frac{d}{d\omega}\right)\left(\frac{f(\omega)}{\omega}\right)^s \cdot \sum_{n=0}^\infty \frac{g_n \omega(x)^n}{n!}=\sum_{n=0}^\infty\frac{s-n-1}{n+1}\left(g_{n+1}-g_0q_{n+1}(s)\right)\frac{\omega(x)^n}{n!}
$$
\emph{Proof}: Suppose $n>0$. Then
\begin{align*}
&\left(\frac{\omega}{f(\omega)}\right)^s\left(s\mathrm{L}-\frac{d}{d\omega}\right)\left(\frac{f(\omega)}{\omega}\right)^s \cdot ~\omega(x)^n=\\
&~~~=\left(\frac{\omega}{f(\omega)}\right)^s\left(\frac{s}{x}\omega^n \left(\frac{f(\omega)}{\omega}\right)^s -n\omega^{n-1}\left(\frac{f(\omega)}{\omega}\right)^s -s\omega^n\left(\frac{1}{x}-\frac{1}{\omega}\right)\left(\frac{f(\omega)}{\omega}\right)^s \right)=(s-n)\omega^{n-1}
\end{align*}
Consider the case $n=0$. Then
\begin{align*}
&\left(\frac{\omega}{f(\omega)}\right)^s\left(s\mathrm{L}-\frac{d}{d\omega}\right)\left(\frac{f(\omega)}{\omega}\right)^s \cdot ~\omega(x)^n=\left(\frac{\omega}{f(\omega)}\right)^s\left(s\mathrm{L}-\frac{d}{d\omega}\right)\left(\frac{f(\omega)}{\omega}\right)^s \cdot ~1=\\
&~~~=\left(\frac{\omega}{f(\omega)}\right)^s\left(\frac{s}{x}\left(\frac{f(\omega)}{\omega}\right)^s -\frac{s}{x}-s\left(\frac{1}{x}-\frac{1}{\omega}\right)\left(\frac{f(\omega)}{\omega}\right)^s \right)=\frac{s}{\omega}-\frac{s f'(\omega)}{f(\omega)}\left(\frac{\omega}{f(\omega)}\right)^s
\end{align*}
But
$$
\sum_{n=0}^\infty \frac{q_n(s)\omega^n}{n!}=\left(\frac{\omega}{f(\omega)}\right)^s ~\Rightarrow ~\frac{s}{\omega}-\frac{s f'(\omega)}{f(\omega)}\left(\frac{\omega}{f(\omega)}\right)^s=\sum_{n=0}^\infty \frac{(n+1-s)q_{n+1}(s)\omega^n}{(n+1)!}
$$
And hence
\begin{align*}
&\left(\frac{\omega}{f(\omega)}\right)^s\left(s\mathrm{L}-\frac{d}{d\omega}\right)\left(\frac{f(\omega)}{\omega}\right)^s \sum_{n=0}^\infty\frac{g_n}{n!}\omega^n=\\
&~~~=\sum_{n=1}^\infty \frac{(s-n)g_n}{n!}\omega^{n-1}+g_0\sum_{n=0}^\infty \frac{(n+1-s)q_{n+1}(s)\omega^n}{(n+1)!}=\\
\tag*{\qed}&~~~=\sum_{n=0}^\infty\frac{s-n-1}{n+1}\left(g_{n+1}-g_0q_{n+1}(s)\right)\frac{\omega(x)^n}{n!}
\end{align*}
\textbf{Proposition A.2}. Suppose
$$
\left(\frac{\omega}{f(\omega)}\right)^s\left(s\mathrm{L}-\frac{d}{d\omega}\right)^k\left(\frac{f(\omega)}{\omega}\right)^s \cdot ~\sum_{n=0}^\infty \frac{g_n^0 \omega(x)^n}{n!}=\sum_{n=0}^\infty\frac{g_n^k \omega(x)^n}{n!}
$$
Then $\forall n\geqslant 0$ the following holds:
$$
\binom{s-1}{n}(g_n^k-g_0^kq_n(s))+\sum_{m=0}^k\binom{s-1}{n+m}q_{n+m}(s)g_0^{k-m}=\binom{s-1}{n+k}g_{n+k}^0
$$
\emph{Proof}:
\begin{align*}
&g_n^k=\frac{s-n-1}{n+1}g_{n+1}^{k-1}-\frac{s-n-1}{n+1}g_0^{k-1}q_{n+1}(s) ~\Rightarrow\\
&~~~\Rightarrow~g_n^k+\frac{s-n-1}{n+1}g_0^{k-1}q_{n+1}(s)=\frac{s-n-1}{n+1}\left(\frac{s-n-2}{n+2}g_{n+2}^{k-2}-\frac{s-n-2}{n+2}g_0^{k-2}q_{n+2}(s)\right)
\end{align*}

\newpage
\begin{center}
\textbf{Appendix A}
\end{center}

\hrule
~\\
~\\
Repeatedly using this argument we obtain (here $(x)_n$ denotes the falling factorials):
\begin{align*}
g_n^k+\frac{(s-n-1)_1}{(n+1)_1}g_0^{k-1}q_{n+1}(s)&+\frac{(s-n-1)_2}{(n+2)_2}g_0^{k-2}q_{n+2}(s)+...+\frac{(s-n-1)_k}{(n+k)_k}g_0^{0}q_{n+k}(s)=\\
=\frac{(s-n-1)_k}{(n+k)_k}g_{n+k}^{0}
\end{align*}
Now multiply by $\binom{s-1}{n}$ to rewrite the latter as
$$
\binom{s-1}{n}g_{n}^{k}+\binom{s-1}{n+1}g_0^{k-1}q_{n+1}(s)+...+\binom{s-1}{n+k}g_0^0q_{n+k}(s)=\binom{s-1}{n+k}g_{n+k}^0
$$
And that proves the statement. Particularly, in case $n=0$, since $q_0(s)=1$, we have:
\begin{align*}
\tag*{\qed}\sum_{m=0}^k\binom{s-1}{m}q_{m}(s)g_0^{k-m}=\binom{s-1}{k}g_{k}^0
\end{align*}
It then follows that
$$
\sum_{k=0}^\infty g_0^{k}\alpha^{-k}~\cdot~\sum_{k=0}^\infty\binom{s-1}{k}q_k(s)\alpha^{-k}=\sum_{k=0}^\infty \binom{s-1}{k}g_k^0\alpha^{-k}
$$
Hence,
$$
\sum_{k=0}^\infty \binom{s-1}{k}q_k(s)\alpha^{s-k}=p_s(\alpha) ~~\Rightarrow ~~\sum_{k=0}^\infty g_0^k\alpha^{-k}=\frac{\alpha}{p_s(\alpha)}\sum_{k=0}^\infty\frac{g_k^0}{k!}D^k \cdot \alpha^{s-1}
$$
Let $\ell_s(\alpha)$ denote the resulting series $\sum_{k\geqslant0} g_0^k\alpha^{-k}$. Then for an arbitrary $n$ the following holds:
$$
\sum_{k=0}^\infty g_n^k\alpha^{-k}=\frac{1}{\binom{s-1}{n}}\left[q_n(s)\binom{s-1}{n}\ell_s(\alpha)-\ell_s(\alpha)\sum_{k=0}^\infty \binom{s-1}{n+k}q_{n+k}(s)\alpha^{-k}+\sum_{k=0}^\infty\binom{s-1}{n+k}g_{n+k}^0 \alpha^{-k}\right]
$$
Thus
\begin{align*}
&\left(\frac{\omega}{f(\omega)}\right)^s\sum_{k=0}^\infty\alpha^{-k}\left(s\mathrm{L}-\frac{d}{d\omega}\right)^k\left(\frac{f(\omega)}{\omega}\right)^s \cdot\sum_{n=0}^\infty\frac{g_n^0}{n!}\omega^n=\sum_{k=0}^\infty\alpha^{-k}\sum_{n=0}^\infty\frac{g_n^k}{n!}\omega^n=\\
&~~~~~~~=\sum_{n=0}^\infty\frac{\omega^n}{n!\binom{s-1}{n}}\left[\ell_s(\alpha)q_n(s)\binom{s-1}{n}+\sum_{k=0}^\infty\binom{s-1}{n+k}(g_{n+k}^0-\ell_s(\alpha)q_{n+k}(s)) \alpha^{-k}\right]=\\
&~~~~~~~=\left(\frac{\omega}{f(\omega)}\right)^s\ell_s(\alpha)+\sum_{n=0}^\infty\frac{\omega^n}{(s-1)_n}\sum_{k=0}^\infty\binom{s-1}{n+k}(g_{n+k}^0-\ell_s(\alpha)q_{n+k}(s)) \alpha^{-k}
\end{align*}
The latter may be rewritten as
\begin{align*}
&\left(\frac{\omega}{f(\omega)}\right)^s\ell_s(\alpha)+\underline{\sum_{n=0}^\infty\frac{(\alpha\omega)^n}{(s-1)_n}\sum_{k=0}^\infty\binom{s-1}{k}(g_k^0-\ell_s(\alpha)q_k(s))\alpha^{-k}}-\\
&~~~~~~-\sum_{n=1}^\infty\frac{(\alpha\omega)^n}{(s-1)_n}\sum_{k=0}^{n-1}\binom{s-1}{k}(g_k^0-\ell_s(\alpha)q_k(s))\alpha^{-k}
\end{align*}
By definition of $\ell_s(\alpha)$, the second summand is equal to zero and hence this expression may be rewritten again as
$$
\left(\frac{\omega}{f(\omega)}\right)^s\ell_s(\alpha)-\sum_{k=0}^\infty\frac{(g_k^0-\ell_s(\alpha)q_k(s))}{k!}\alpha^{-k}\sum_{n=k}^\infty\frac{(\alpha\omega)^{n+1}(s-1)_{k}}{(s-1)_{n+1}}
$$

\newpage
\begin{center}
\textbf{Appendix A}
\end{center}

\hrule
~\\
~\\
Now notice that the following formally holds (in fact only for $s$: $\mathrm{Re}~s < 1+k$):
$$
\sum_{n=k}^\infty \frac{x^n (s-1)_k}{(s-1)_{n+1}}=x^k\sum_{n=0}^\infty \frac{(-1)^{n+1}x^n}{n!}\frac{\Gamma(1+n)\Gamma(1+k-s)}{\Gamma(2+n+k-s)}=-x^k\int_{0}^{1}t^{k-s}e^{-x(1-t)}dt
$$
It then follows that for the series $\sum_{n\geqslant 0}g_n^0 x^n/n!=g(x)$ holds:
\begin{align*}
&\left(\frac{\omega}{f(\omega)}\right)^s\sum_{k=0}^\infty\alpha^{-k}\left(s\mathrm{L}-\frac{d}{d\omega}\right)^k\left(\frac{f(\omega)}{\omega}\right)^s \cdot~g(\omega(x))=\\
&~~~~~~~=\left(\frac{\omega}{f(\omega)}\right)^s\ell_s(\alpha)+\alpha\omega(x)\sum_{k=0}^\infty\frac{(g_k^0-\ell_s(\alpha)q_k(s))}{k!}\omega^{k}\int_{0}^{1}t^{k-s}e^{-\alpha\omega(x)(1-t)}dt=\\
&~~~~~~~=\left(\frac{\omega}{f(\omega)}\right)^s\ell_s(\alpha)+\alpha\omega(x)\int_{0}^{1}t^{-s}e^{-\alpha\omega(x)(1-t)}\left(g(t\omega(x))-\ell_s(\alpha)\left(\frac{t\omega(x)}{f(t\omega(x))}\right)^s\right)dt
\end{align*}
By definition,
$$
\ell_s(\alpha)=\dfrac{\alpha}{p_s(\alpha)}g(D)\cdot\alpha^{s-1}=\dfrac{\alpha}{p_s(\alpha)}g(D)\left(\dfrac{f(D)}{D}\right)^s\cdot\dfrac{p_s(\alpha)}{\alpha}
$$
Hence one may multiply and divide by $\left(f(\omega)/\omega\right)^s$, to obtain
\begin{align*}
&\sum_{k=0}^\infty\alpha^{-k}\left(s\mathrm{L}-\frac{d}{d\omega}\right)^k\cdot~ g(\omega(x))=\dfrac{\alpha}{p_s(\alpha)}g(D)\cdot\frac{p_s(\alpha)}{\alpha}+\\
&~~~~~~~~+\alpha\omega(x)e^{-\alpha\omega(x)}\int_{0}^{1}e^{\alpha t\omega(x)}\left(\frac{f(\omega(x))}{f(t\omega(x))}\right)^s\left(g(t\omega(x))-\dfrac{\alpha}{p_s(\alpha)}g(D)\cdot\frac{p_s(\alpha)}{\alpha}\right)dt
\end{align*}
or equivalently
\begin{align*}
&\left(1-s\alpha^{-1}\mathrm{L}+\alpha^{-1}\frac{d}{d\omega}\right)^{-1} \cdot~ g(\omega(x))=\\
&~~~~~=\frac{\alpha}{p_s(\alpha)}g(D)\cdot\frac{p_s(\alpha)}{\alpha}+\alpha e^{-\alpha\omega(x)}f(\omega(x))^s\int_{0}^{\omega(x)}e^{\alpha t}f(t)^{-s}\left(g(t)-\frac{\alpha}{p_s(\alpha)}g(D)\cdot\frac{p_s(\alpha)}{\alpha}\right)dt
\end{align*}
The latter expression may also be obtained by the use of hypergeometric-like representation (1.11), since one may solve the differential equation of the second order, and separately find the value of the solution at zero.\\
\phantom{q}\hfill\ensuremath{\blacksquare}

\end{document}